\documentclass[11pt]{amsart}
\usepackage{amssymb}
\usepackage{graphicx}

\newcommand{\Aut}{{\mathrm{Aut} }}

\newcommand{\ad}{{\mathrm{ad} }}
\newcommand{\F}{\mathbf F}
\newcommand{\vol}{{\mathrm{vol}\, }}
\newcommand{\R}{\mathbf R}

\renewcommand{\deg }{{\mathrm{deg\,} }}
\renewcommand{\max }{{\mathrm{max\,} }}

\newcommand{\GL}{{\mathbf{GL} }}

\newcommand{\rad}{{\mathrm{rad}\, }}
\newcommand{\SL}{{\mathbf{SL} }}

\newcommand{\Tr}{{\mathrm{Tr}\, }}
\newcommand{\Z}{\mathbf Z}
\newcommand{\C}{\mathbf C}
\newcommand{\spa}{\mathrm{span}\,}
\renewcommand{\b}[1]{{\mathbf{#1}}}

\renewcommand{\deg}{{\mathrm{deg\,} }}

\renewcommand{\b}[1]{{\mathbf{#1}}}

\newcommand{\h}{\mathbf{F}_q}

\newcommand{\pp}{\mathfrak{p}}

\renewcommand{\text}{\mathrm}

\newtheorem{theorem}{Theorem}[section] 
\newtheorem{lemma}[theorem]{Lemma}
\newtheorem{corollary}[theorem]{Corollary} 
\newtheorem{proposition}[theorem]{Proposition}

\theoremstyle{definition}
\newtheorem{definition}{Definition}

\begin{document}

\title{Isospectral Definite Ternary $\h[t]$-Lattices}
\author{Jean Bureau}
\address{Mathematics Department, Louisiana State University, Baton Rouge, LA 70803-4918, U.S.A.}
\email{jbureau@math.lsu.edu} 
\author{Jorge Morales}
\address{Mathematics Department, Louisiana State University, Baton Rouge, LA 70803-4918, U.S.A.}
\email{morales@math.lsu.edu} 

\begin{abstract}
We prove that the representations numbers of a ternary definite integral
quadratic form defined over $\h[t]$, where $\h$ is a finite field of odd
characteristic, determine its integral equivalence class when $q$ is
large enough with respect to its successive minima. Equivalently, such a
quadratic form is determined up to integral isometry by its theta series.
\end{abstract}

\maketitle

\section{Introduction}

It has been a long-standing question to determine whether integral definite $\Z $-lattices are determined up
to isometry by their theta series. In 1979, Watson \cite{Watson:1979we} proved that definite binary $\Z $-lattices are determined by their primitive representations. The case of ternary lattices had to wait until 1997 to be solved by Schiemann \cite{Schiemann:1997rf} by means of extensive computations. He proved that definite ternary $\Z $-lattices are indeed determined by their representation numbers. This is not the case for forms of rank $\ge 4$, where counterexamples have been found (cf. \cite{Schiemann:1990zn}, \cite{Kitaoka:1977or} and \cite{Conway:1992vo}). 

In this paper we prove the analogue of Schiemann's theorem for definite ternary $\h[t]$-lattices, where $\h$ is a finite field of odd characteristic. We show first that the representation
numbers determine invariants such as the successive minima and the genus (Sections \ref{S:minima} and \ref{S:genus}). 
Our proof that the representation
numbers determine the equivalence class requires different arguments according to different configurations
of the successive minima (Section \ref{S:equivalence}). When the successive minima have alternating parity, we use a theorem of Carlitz based on Fourier inversion and we are able to conclude
equivalence under the hypothesis that the ground field $\h $ is large enough (see Theorem \ref{T:case3} for a precise statement). This condition is not required
in the two other cases (Theorems \ref{T:case1} and \ref{T:case2}).

We thank the referee for his/her careful reading and useful comments.

\section{Notation and terminology}

The following notation will be in force throughout the paper:

\begin{itemize}
\item[$\h$ :] The finite field of order $q$. We always assume $q$ odd.
\item[$A$ :] The polynomial ring $\h[t]$
\item[$K$ :] The field of rational functions $\h(t)$
\item[$\delta$ :] A fixed non-square of $\h^\times$.
\end{itemize}

Let $L$ be an $A$-lattice of finite rank $n$ and let $Q$ be
a quadratic form on $L$. The form $Q$ is {\em definite} if it is anisotropic over
the field $K_\infty=\h((1/t))$. This implies in particular
that $n\le 4$. 

Let $B(\b x,\b y)=Q(\b x+\b y)-Q(\b x)-Q(\b y)$ be the associated symmetric bilinear form. 
Djokovi{\'c} \cite{Djokovic:1976ve} showed that if $(L,Q)$ is
definite, then there exists an $A$-basis ${\b v}_1,\ldots,{\b v}_n$ of $L$
such that the Gram matrix $M=(m_{ij})$, where $m_{ij}=\frac{1}{2}B({\b v}_i,{\b v}_j)$, satisfies 
\begin{equation}\label{E:reduced}
\deg m_{ii}\le \deg
m_{jj}\ \mathrm{for}\  i\le j \ \mathrm{and}\ \deg m_{ij}< \deg m_{ii}\ \mathrm{for}\  i<j.
\end{equation}
 
Such a basis is called {\em reduced}. Gerstein \cite[Theorem 2]{Gerstein:2003rf} showed
that if ${\b v}_1',\ldots,{\b v}_n'$ is another reduced basis for $(L,Q)$, then
\[
{\b v}_j'=\sum_{i=1}^n u_{ij} {\b v}_i,
\]
where $u_{ij}\in \h$.

In
particular, the increasing sequence of degrees of the diagonal terms of
a reduced Gram matrix
\[
(\deg m_{11}, \deg m_{22},\ldots,\deg m_{nn})
\]
is an invariant of the equivalence class of the quadratic form. This sequence
is called the {\em sequence of successive minima} of $Q$ and will be denoted
by \[
(\mu_1(L,Q),\mu_2(L,Q),\ldots,\mu_n(L,Q)).
\]

The {\em representation numbers} of $(L,Q)$ are defined by
\begin{equation}\label{E:reprno}
R(L,Q,a)=|\{{\b v} \in L: Q({\b v} )=a\}|\quad (a\in K).
\end{equation} 

It is easy to see that if $(L,Q)$ is definite, the above numbers
are finite. Clearly they depend only on the isometry class
 of $(L, Q)$. 
\\
 
\begin{definition} Two definite quadratic forms $(L,Q)$ and $(L',Q')$
are called {\em isospectral}\footnote{The terminology 
comes from the fact that for quadratic forms over $\Z$
the representation numbers are naturally the dimensions of
the eigenspaces of a Laplace operator, see \cite{Milnor:1964ao}.} if $R(L,Q,a)=R(L',Q',a)$ for all $a\in K$.
\end{definition}

Following Conway's \cite{Conway:1997hc} terminology, we shall call an invariant
of $(L,Q)$ {\em audible} if it is determined by the representation
numbers. The main goal of this paper is to show that the equivalence
class of a ternary definite quadratic form over $A$ is audible. We
shall do this in several steps.

\section{The successive minima}\label{S:minima}

Let $(L,Q)$ be a definite quadratic form over $A$. For $m\in\Z$,
define
\begin{equation}\label{E:filter}
L_m=\{{\b v}  \in L: \deg Q({\b v} )\le m\}.
\end{equation}
It is easy to see that the $L_m$ are finite-dimensional
$\h$-subspaces of $L$ and that they form an increasing sequence
whose union is $L$. We encode their successive dimensions into the
formal power series
\begin{equation}\label{E:pseries}
{\b S}_L(u)=\sum_{m\in \Z} \dim (L_m/L_{m-1}) u^m\quad\mathrm{and}\quad {\b T}_L(u)= \sum_{m\in \Z} \dim (L_m) u^m.
\end{equation}
Notice that both ${\b S}_L(u)$ and ${\b T}_L(u)$ are Laurent series in $u$ since 
$L_m=\{0\}$ for $m\ll 0$ (we do not assume that $Q$ takes integral values on $L$). It is clear from their definition that both series are audible.

\begin{proposition}\label{P:minima} With the notation above we have
\[
{\b S}_L(u)=\frac{u^{\mu_1}+u^{\mu_2}+u^{\mu_3}}{1-u^2}\quad\mathrm{and}\quad
{\b T}_L(u)=\frac{u^{\mu_1}+u^{\mu_2}+u^{\mu_3}}{(1-u^2)(1-u)},
\]
where $(\mu_1,\mu_2,\mu_3)$ are the successive minima
of $(L,Q)$. In particular the sequence $(\mu_1,\mu_2,\mu_3)$ is audible.
\end{proposition}

\begin{proof} Let ${\b v}_2, {\b v}_2, {\b v}_3$ be a reduced basis of $L$.
Notice that since $Q$ is definite, $\mu_1,\mu_2, \mu_3$ cannot
all have the same parity. 

Suppose first that $\mu_1\equiv \mu_2 \pmod 2$. If $m<\mu_1$, then clearly
$L_m$ is trivial. When $m\equiv \mu_1 \pmod 2$ and $\mu_1\le m < \mu_2$, the quotient space
$L_m/L_{m-1}$ is 1-dimensional (with basis $\{t^{(m-\mu_1)/2} \b
v_1\}$). When $m\equiv \mu_1 \pmod 2$ and $m\ge \mu_2$,
the quotient $L_m/L_{m-1}$ is $2$-dimensional (with basis $\{t^{(m-\mu_1)/2} \b
v_1 , t^{(m-\mu_2)/2} {\b v}_2\}$). When $m\equiv \mu_3 \pmod 2$, the quotient $L_m/L_{m-1}$  is trivial if $m<\mu_3$
and $1$-dimensional if $m\ge \mu_3$ (with basis $\{t^{(m-\mu_3)/2}
{\b v}_3\}$).

Putting this information into the series, we get

\[
\begin{aligned}
{\b S}_L(u)= & \sum_{k=0}^{(\mu_2-\mu_1)/2-1} u^{\mu_1+2 k} + 2\sum_{k=0}^{\infty} u^{\mu_2+2 k}
+\sum_{k=0}^{\infty} u^{\mu_3+2 k}\\
= & \frac{u^{\mu_1}+u^{\mu_2}+u^{\mu_3}}{1-u^2}.
\end{aligned}
\]

The case when $\mu_1\not\equiv \mu_2 \pmod 2$ is computed similarly. We spare the reader the details. The second identity follows from 
the obvious relation ${\b S}_L(u)=(1-u){\b T}_L(u)$.

\end{proof}

\section{The genus}\label{S:genus}

Let $\pp$ be a prime ideal of $A$ and let $\xi$ be a root (in an algebraic closure of $\F_q$) 
of a generator of $\pp$.  The {\em canonical character} $\chi_{\pp}:
K_{\pp}\rightarrow \mathbb{C}^{\times}$ is the homomorphism  defined by
\[
\chi_{\pp} (f)= \exp\left(2\pi i \Tr (\mathrm{Res}_\xi(f))/p\right),
\]
where $\mathrm{Res}_\xi(f)\in \F_q(\xi)$ is the residue of $f$ at $\xi$ (i.e. the coefficient of $(T-\xi)^{-1}$
in the Laurent series expansion of $f$ at $\xi$) and $\Tr: \F_q(\xi)\to \F_p$ is the trace to the prime field $\F_p$. Clearly the definition is independent of the choice
of the root $\xi$, since residues at different roots are conjugate over $\F_q$. Notice that $\chi_\pp$ is trivial on $A_\pp$; in fact $A_\pp$ is the largest 
fractional ideal of $K_\pp$ on which $\chi_\pp$ is trivial.\\

Let $(W,Q)$ be a definite quadratic space over $K$ and let $L\subset W$ be an $A$-lattice,
not necessarily integral with respect to $Q$.

Define
\[
\mu(L,Q,\chi_{\pp})=\lim_{m\rightarrow \infty}\frac{1}{|L_m|}\>\sum_{\b x\in L_m}\chi_{\pp}(Q(\b x)).
\]
We shall see below that this is a stabilizing limit. We first notice that this ``average'', $\mu(L,Q,\chi_{\pp})$, is audible. Indeed, we have 
\begin{align*}
\label{mu}
\mu(L,Q,\chi_{\pp})&=\lim_{m\rightarrow \infty}\frac{1}{|L_m|}\>\sum_{\b x\in L_m}\chi_{\pp}(Q(\b x))\\
&=\lim_{m\rightarrow \infty} \frac{1}{|L_m|}\sum_{\deg (a)\leq m}R(L,a)\chi_{\pp}(a).\\
\end{align*}
We now express $\mu(L,Q,\chi_{\pp})$ in terms of local data. Let $L^\sharp$ 
be the dual of $L$ with respect to $Q$. Since $L$ is the union of the $L_m$, 
for $m$ large enough, the restriction $ L_m\to L/(L^\sharp \cap L)$ of the canonical
projection  is surjective. Thus

\begin{equation}\label{E:avg}
\begin{aligned}
\frac{1}{|L_m|}\>\sum_{\b x\in L_m}\chi_{\pp}(Q(\b x))& = \frac{|L_m\cap L^\sharp|}{|L_m|}
\sum_{\b x\in L_m/L_m\cap L^\sharp}\chi_{\pp}(Q(\b x))\\
&= \frac{1}{|L:L\cap L^\sharp|}\sum_{\b x\in L/L\cap L^\sharp}\chi_{\pp}(Q(\b x))\\
&= \frac{1}{|L_\pp:L_\pp\cap L_\pp^\sharp|}\sum_{\b x\in L_\pp/L_\pp\cap L^\sharp_\pp}\chi_{\pp}(Q(\b x)).
\end{aligned}
\end{equation}

\begin{theorem}\label{T:audible}
Let $\pi$ be a monic generator of $\pp$. The
sequence $\mu(L,\pi^{-k} Q,\chi_{\pp})$ $(k=0,1,2,\ldots)$ 
determines completely the local class $(L_\pp,Q)$. 
\end{theorem}
\begin{proof}
Let $(L_\pp,Q)=(M_1,Q_1)\perp (M_2,Q_2)\perp\ldots\perp (M_r,Q_r)$ be the Jordan decomposition. 
Each 
$M_i$ is $\pp$-modular, i.e. $M_i^\sharp=\pi^{-\nu_i} M_i$, and we assume $\nu_1<\nu_2<\cdots <\nu_r$.
We define $\mu$ for local lattices using the last line of \eqref{E:avg}.
Then we have
\[
\mu(M_i,\pi^{-k}Q_i,\chi_{\pp})=\begin{cases} 
\displaystyle{[M_i:\pi^{k-\nu_i}M_i]^{-1}\sum_{\b x\in M_i/\pi^{k-\nu_i}M_i}\chi_p(\pi^{-k}Q(\b x))}& \mathrm{if }\quad  k\ge \nu_i ; \\
\qquad 1 & \mathrm{if }\quad  k< \nu_i.
\end{cases}
\]
We can express this further using the 
invariant $\gamma_\pp$ defined in \cite[Ch. V, \S8]{Scharlau:1985jv} (see also \cite[
\S 24]{Weil:1964ay}). Then we have

\begin{equation}\label{E:weil}
\mu(M_i,\pi^{-k}Q_i,\chi_{\pp})=
\begin{cases} 
[M_i:\pi^{k-\nu_i}M_i]^{-1/2}\gamma_\pp(\pi^{-k} Q_i)& \mathrm{if }\quad  k\ge \nu_i ; \\
\qquad 1 & \mathrm{if }\quad  k< \nu_i.
\end{cases}
\end{equation}

The invariant $\gamma_\pp(\pi^{-k} Q_i)$ is a 4th root of unity and depends only
on the class of $\pi^{-k} Q_i$ over the field $K_\pp$ (actually, only on its Witt class)
\cite[Chapter 5]{Scharlau:1985jv}. In particular $|\gamma_\pp(\pi^{-k} Q_i)|=1$, thus
\[
\log_q|\mu(M_i,\pi^{-k}Q_i,\chi_{\pp})|= - \frac{m_i\deg \pi}{2} \sup\{0, k-\nu_i\},
\]
where $m_i$ is the rank of $M_i$ and $\log_q$ is the logarithm in base $q$. 
Using the obvious fact that $\mu$ is compatible
with orthogonal sums, we get

\begin{equation}\label{E:eqn0}
\log_q|\mu(L,\pi^{-k} Q,\chi_{\pp})|=- \sum_{i=1}^r \frac{m_i\deg \pi}{2} \sup\{0, k-\nu_i\}.
\end{equation}
As observed earlier, the left-hand side of \eqref{E:eqn0} is audible as a function of $k$, then so is the 
right-hand side. The functions $f_\nu: \Z \to \R$ given by $f_\nu(k)=\sup\{0,k-\nu\}$ are linearly independent,
so the expression of $\log_q|\mu(L,\pi^{-k} Q,\chi_{\pp})|$ in \eqref{E:eqn0} as linear combination of these functions 
is unique; it follows that
 the numbers $\nu_1,\nu_2,\ldots,\nu_r$ and the ranks $m_1,m_2,\ldots,m_r$ of the Jordan factors
of $L_\pp$ are audible.

It is left to show that $\det (M_i,Q_i)$ is audible. Consider the case $i=1$ and $k=\nu_1+1$. Let
$F=\pi^{-\nu_1} Q_1$ (note that $F$ is unimodular on $M_1$). By
\eqref{E:weil}, we have
\begin{equation}\label{E:weil2}
\begin{aligned}
\mu(L,\pi^{-\nu_1-1}Q,\chi_{\pp})&=\mu(M_1,\pi^{-1} F,\chi_{\pp})\\
&= q^{-m_1 \deg\pi/2}\gamma_\pp(\pi^{-1} F).
\end{aligned}
\end{equation}
The invariant $\gamma_\pp$ satisfies $\gamma_\pp(<a>)\gamma_\pp(<b>)=\gamma_\pp(<ab> )(a,b)_\pp$,
 where $(a,b)_\pp$ is the Hilbert symbol \cite[\S28 p. 176]{Weil:1964ay}. Applying this identity, we get
\[
\gamma_\pp(\pi^{-1} F)=\gamma_\pp(\langle\pi\rangle)^{m_1}\gamma_\pp(F) (\det F,\pi)_\pp,
\]
where $\langle\pi\rangle$ is the rank-one form $\pi X^2$. Since $F$ is unimodular on $M_1$,
$\gamma_\pp(F) =1$, so $\gamma_\pp(\pi^{-1} F)=\gamma_\pp(\langle\pi\rangle)^{m_1} (\det F,\pi)_\pp$.
It follows from this and \eqref{E:weil2} that 
the class of $(M_1,Q_1)$ is audible. We continue similarly taking successively 
$k=\nu_2+1,\ldots, \nu_r+1$.

\end{proof}

Theorem \ref{T:audible} has two immediate consequences:

\begin{corollary}\label{C:genus} The genus of $(L,Q)$ is audible.
\end{corollary}
\begin{corollary}\label{C:disc} The discriminant of $(L,Q)$ is audible.
\end{corollary}

\section{The theta series and the adjoint form}\label{S:theta}

Let $(L,Q)$ be a definite ternary $A$-lattice. We define the theta series 
of $(L,Q)$ as in R\"uck \cite{Ruck:1995dt} and Rosson \cite{Rosson:2002tn}. 
We shall refer to these
papers for details of some computations. 

The analogue of the Poincar\'e complex half-plane is $\mathfrak{H}=\SL_2(K_{\infty})/\SL_2(\mathcal{O}_{\infty})$.
A complete set of coset representatives for $\mathfrak{H}$ is the set 
\begin{equation}\label{E:cosets}
\mathfrak{D}=\left\{\left [\begin{matrix} y&xy^{-1}\\0&y^{-1}\end{matrix}\right ] : y=t^m,\  m\in \Z,\ x\in t^{2m+1} A\right\}.
\end{equation}
Let $x=\sum_{i=-\infty}^n x_i t^{i}\in K_{\infty}$. We define a character of $e:K_{\infty}\to \C^\times$ by 
\[
e\{x\}=\exp(2i\pi\,\Tr(x_1)/p),
\]
where $\Tr$ stands for the trace of $\h$ to its prime subfield and $p$ is the characteristic of $\h$.
 Let $\Psi$ denote the characteristic function of $\mathcal{O}_{\infty}$. For $z=\left(\begin{matrix}y&xy^{-1}\\0&y^{-1}\end{matrix}\right)\in \mathfrak{H}$ and for a lattice $L$, we define the theta series of $L$ by
\begin{align*}
\theta_L(z)&=\sum_{{\b v} \in tL}\Psi\left(y^2Q({\b v} )\right)\>e\{xQ({\b v} )\}\\
&=\sum_{w\in L}\Psi(t^2y^2Q(\b w))\>e\{t^2xQ(\b w)\}.
\end{align*} 

It is readily checked that $\theta_L$ is a function on $\mathfrak{H}$, i.e. does not depend
on the chosen matrix representatives. The theta series determines the representation numbers and conversely. 
Indeed, for $y=t^{-m}$, we have  
\begin{equation}
\label{E:use}
\theta_L(z)=\sum_{v_{\infty}(a)\geq 2m-2}R(L,a)\>e\{xt^2 a\}.
\end{equation}
It is clear from this that the representation numbers $R(L,a)$ can be
recovered from $\theta_L(z)$ by Fourier inversion. 

Let $d{\b v} $ be an
additive Haar measure on $V_{\infty}$. For a locally constant compactly
supported function $f$ on $V_{\infty}$, we define its Fourier 
transform by 
\[
\hat{f}(\b w)=\int_{V_{\infty}} f({\b v} )e\{-B({\b v} ,\b w)\}d{\b v} ,
\]
where $B$ is the bilinear form associated to $Q$. We shall further
assume that the Haar measure $d {\b v} $ is self-dual, i.e. it has been normalized 
so that
\begin{equation}\label{E:normalization}
\hat{\hat{f}}({\b v} )=f(-{\b v} ).
\end{equation}
This is equivalent to saying that the volume with respect to $d{\b v} $ of any $O_\infty$-lattice $M\subset V_\infty$ 
satisfies
\[
\vol(M)\vol(M^*)=1,
\]
where $M^*=\{\b w\in V_\infty: B(\b w, M)\subset O_\infty\}$.

\begin{proposition}\label{P:transform}
Let $G,H\in K_{\infty}$ , $H\ne 0$,  be such that $v_{\infty}(G)=g>h=v_{\infty}(H)$. 
Let $\varphi: V_\infty \to \C$ be the function defined by $\varphi({\b v} )=\Psi(Q({\b v} ) G)e(Q({\b v} )H).$
Then the Fourier transform of $\varphi$ is given by
\begin{equation}\label{E:fourier}
\hat\varphi (\b w) = I \Psi\left(\frac {G}{H^{2}} Q(\b w)\right) e\left(-\frac{1}{H}Q(\b w)\right),
\end{equation}
where $I=  |H|_{\infty}^{-3/2} \gamma_{\infty}(H Q)$.
\end{proposition}

\begin{proof} Essentially the same computation as in \cite[Theorem 3.2]{Rosson:2002tn},
shows \eqref{E:fourier} with
\[
I=\int_{V_{\infty}} \Psi(Q({\b v} ) G)e(Q({\b v} )H) d{\b v} .
\]
We shall evaluate $I$ explicitly. Since $Q$ is definite, there exists
a unique $O_{\infty}$-lattice $M\subset V_{\infty}$ maximal
with respect to the property $G Q(M)\subset O_{\infty}$. Then
\[
I=\int_{M} e(Q({\b v} )H) d{\b v} .
\]
We shall now see that the form $HQ$ is integral on $H^{-1} M^{*}$.
On the one hand, since $g>h$ we have $H^{-1} M^{*}=(H^{-1}G) (G^{-1}M^{*})\subset t^{-1} G^{-1}M^{*}$. On the other hand, since $M$ is maximal integral with respect to $G Q$, we have $t^{-1} G^{-1}M^{*}\subset M$. Thus
\[
I=\vol(H^{-1} M^{*})[M:H^{-1} M^{*}]^{1/2}\gamma_{\infty}(HQ).
\]
To finish the computation, we observe
\[
\begin{aligned}
\vol(H^{-1} M^{*})[M:H^{-1} M^{*}]^{1/2}&= \vol(H^{-1} M^{*})^{1/2}\vol(M)^{1/2}\\
&=|H|_{\infty}^{-3/2}[\vol(M^{*})\vol(M)]^{1/2}\\
&=|H|_{\infty}^{-3/2}.
\end{aligned}
\]
Notice that the last line uses the chosen normalization
\eqref{E:normalization} for the Haar measure.

\end{proof}

\begin{corollary}\label{C:poisson} Let $z=\left [\begin{matrix} y&xy^{-1}\\0&y^{-1}\end{matrix}\right ] \in \mathfrak{D}$ with $x\ne 0$ and let $S=\left [\begin{matrix} 0&-1\\1&0\end{matrix}\right ] $. Then

\[
\theta_{L}(z)= |D|^{-1/2}_{{\infty}} I(z) \theta_{L^{\sharp}}(S\cdot z),
\]
where $I(z)= |x|_{\infty}^{-3/2} \gamma_{\infty}(x Q).$

\end{corollary}

\begin{proof} 
 Let $G=y^2$ and $H=x$. Since $z\in \mathfrak{D}$
we have $v_\infty(y^2)>v_\infty(x)$, so $G$ and $H$ satisfy
the hypotheses of Proposition \ref{P:transform}. Moreover
\[
S\cdot z\sim \left [\begin{matrix} 0&-1\\1&0\end{matrix}\right ] \left [\begin{matrix} y&xy^{-1}\\0&y^{-1}\end{matrix}\right ]  \left [\begin{matrix} 
 1 & 0 \\
 -y^2x^{-1} & 1
\end{matrix}\right ] =\left [\begin{matrix} 
 yx^{-1} & -x^{-1}y^{-1}x \\
 0 & y^{-1}x
\end{matrix}\right ] ,
\]
so applying Proposition \ref{P:transform} to the function
\[
\varphi_z({\b v} )=\Psi(Q({\b v} ) y^2)e(Q({\b v} )x)
\]
we get
\[
\hat\varphi_z({\b v} )=I(z) \varphi_{S\cdot z}({\b v} ).
\]
Applying the Poisson summation formula, we obtain
\[
\sum_{{\b v} \in t L}\varphi_z({\b v} )=\vol(V_\infty/t L)^{-1} \sum_{\b w\in t L^\sharp}\hat\varphi_z({\b v} ),
\]
hence
\
\[
\theta_{L}(z)= |D|^{-1/2}_{{\infty}} I(z) \theta_{L^{\sharp}}(S\cdot z).
\]

(Notice that $\vol(V_\infty/t L)=|D|_\infty ^{1/2}$.)

\end{proof}

Recall that for a ternary lattice $(L,Q)$, its {\em adjoint} $(L^\ad, Q^\ad)$ is defined by
\[
L^\ad= L^\sharp \quad \mathrm{and}\quad Q^\ad=D Q,
\]
where $D=\det(L,Q)$. Alternatively, $(L^\ad, Q^\ad)=(\bigwedge^2 L, \bigwedge^2 Q)$,
where $\bigwedge^2$ is the second exterior power operator.

\begin{theorem}\label{T:adjoint} Let $(L,Q)$ and $(L',Q')$ be isospectral
definite ternary lattices. Then $(L^\ad,Q^\ad)$ and $({L'}^\ad,{Q'}^\ad)$ are isospectral.
\end{theorem}

\begin{proof} Notice that  $R(L^\sharp,Q,a)=R(L^\ad, Q^\ad,D a)$ for all $a\in K$,
so it is enough to prove that $\theta_{L^\sharp}=\theta_{{L'}^\sharp}$.

Since $L$ and $L'$ are in the same genus by Corollary \ref{C:genus},
we have $\det(L,Q)=\det(L',Q')$ and $\gamma_{\infty}(x Q)=\gamma_{\infty}(x Q')$. So, by Corollary \ref{C:poisson}, 
$\theta_{L^\sharp}(z)=\theta_{{L'}^\sharp}(z)$ for $x\ne 0$. 

It remains to prove that $\theta_{L^\sharp}(z)=\theta_{{L'}^\sharp}(z)$ when $x=0$. In this case, letting $y=t^{-m}$ we have, by \eqref{E:use},
\[
\theta_{L^\sharp}(z)=|L^\sharp_{2m-2}|.
\]
These numbers are determined by the series ${\b T}_{L^\sharp}(u)$ defined  in \eqref{E:pseries}, which in turn depends only on the successive minima
of $L^\sharp$ by Proposition \ref{P:minima}. The successive minima of $L^\sharp$ are readily seen to be $(-\mu_3, -\mu_2, -\mu_1)$,
where $(\mu_1,\mu_2,\mu_3)$ are the successive minima of $L$. We conclude by Proposition \ref{P:minima} that $|L^\sharp_{2m-2}|=|{L'}^\sharp_{2m-2}|$.
Thus $\theta_{L^\sharp}(z)=\theta_{{L'}^\sharp}(z)$ for all $z$.
\end{proof}

\section{Equivalence}\label{S:equivalence}

Let $(L,Q)$ and $(L',Q')$ be two isospectral definite ternary lattices over $A$. Our aim in this section is to prove that
they are equivalent.

We already proved  in previous sections that they have the same successive minima $(\mu_1,\mu_2,\mu_3)$ and belong to the same genus -- in particular they have
the same determinant -- and that their adjoints are also isospectral.

\begin{proposition}\label{P:commonbinary} Assume $\mu_3>\mu_2$. Then $L_{\mu_2}$
and $L'_{\mu_2}$ span equivalent binary $A$-lattices.
\end{proposition}

\begin{proof} 
Let $({\b v}_1, {\b v}_2, {\b v}_3)$ and $({\b v}_1', {\b v}_2', {\b v}_3')$  and be reduced bases of $L$ and $L'$ respectively. Let
$M=A{\b v}_1+ A {\b v}_2$ and $M'=A{\b v}_1'+ A {\b v}_2'$. Notice that
$\det(M,Q)$ and is a minimal value for $(L^\ad, Q^\ad)$, and that
it is unique (up to a square in $\h$) with this property since $\mu_1+\mu_2<\mu_1+\mu_3$. 

Since $(L^\ad, Q^\ad)$ and $(L'^\ad, Q'^\ad)$ are isospectral
by Theorem \ref{T:adjoint}, we conclude that $\det(M,Q)=\det(M',Q')$.
Since $\mu_3>\mu_2$, $M_{\mu_2}=L_{\mu_2}=L_{\mu_2}'=M'_{\mu_2}$,.
Thus $M$ and $M'$ represent the same values up to degree
$\mu_2$ and have the same determinant. By \cite[Theorem 4.1]{Bureau:2007qt},
we conclude that $(M,Q)$ and $(M',Q')$ are equivalent.

\end{proof} 

\begin{corollary}\label{C: commonbinary} If $(L,Q)$ and $(L',Q')$ are isospectral lattices
with $\mu_3>\mu_2$, then they have reduced bases such that the 
corresponding reduced Gram matrices have the form

\begin{equation}\label{E:commonbinary}
S=\left [\begin{matrix}  a & b & e \\ b & c& f \\ e & f & g \end{matrix}\right ] 
\quad\mathrm{and}\quad S'=\left [\begin{matrix}  a & b & e' \\ b & c& f' \\ e' & f' & g' \end{matrix}\right ] .
\end{equation}
Furthermore, $g$ and $g'$ may be assumed to have the same leading
coefficient.

\end{corollary}

\begin{proof} Immediate from Proposition \ref{P:commonbinary}.
\end{proof}

\begin{lemma}\label{L:adjoints} Let $S$ and $S'$ be matrices as 
in \eqref{E:commonbinary} representing isospectral forms and assume in addition
$\mu_{1}<\mu_{2}<\mu_{3}$.  Replacing if
necessary the pair $(e,f)$ by $(-e,-f)$ in the matrix $S$, the coefficients 
of $S$ and $S'$ satisfy the relation $a (g-g')= e^2 -e'^2$ and $b (e-e')=a(f-f')$. In particular,
$\deg(g-g')<\deg(e-e')$ and $\deg(f-f')<\deg(e-e').$
\end{lemma}
\begin{proof} The  Gram matrix $S^{\ad}$ of the adjoint $(L^\ad,Q^\ad)$ with respect to the reversed dual basis $\{{\b v}_{3}^{*}, {\b v}_{2}^{*},{\b v}_{1}^{*}\}$ is also reduced
\cite[Lemma 4]{Gerstein:2003rf} and has the form
\[
S^{\ad}=\left [\begin{matrix} 
 a c-b^2 & b e-a f & * \\
 b e-a f & a g-e^2 & * \\
* & * & *
\end{matrix}\right ] .
\]
By Theorem \ref{T:adjoint}, the adjoints $(L^\ad,Q^\ad)$ and $({L'}^\ad,{Q'}^\ad)$ are isospectral, so, by Proposition \ref{P:commonbinary}, the binary lattices
 $M=A {\b v}_{3}^{*}+ A{\b v}_{2}^{*}$ and $M'={A {\b v}_{3}'}^{*}+A\b {v_{2}'}^{*}$
 are equivalent. Since their successive minima are distinct, the
only automorphisms of these lattices are of the form $\mathrm{diag}(\pm1,\pm1)$, so we must have in particular
\begin{equation}
a g-e^2 = a g'-e'^2 \quad\mathrm{and}\quad b e-a f=\pm(b e'-a f').
\end{equation}
Replacing ${\b v}_{3}$ by $-{\b v}_{3}$ if necessary, we can assume 
that the second equality holds with the $+1$ sign. So we get
\begin{equation} \label{E:adj}
a (g-g')= e^2 -e'^2 \quad\mathrm{and}\quad b (e-e')=a(f-f').
\end{equation}
The degree inequalities follow immediately from the fact that $\deg(e+e')<\deg a$ and $\deg b<\deg a$
since $S$ and $S'$  
are reduced.
\end{proof}

\begin{lemma}\label{L:ineq} Let $M=A{\b v}_1+ A {\b v}_2\subset L$. Then
for every $\b w \in M\setminus\{0\}$ we have 
\[
\deg B(\b w,{\b v}_{3})<\deg Q(\b w).
\]
\end{lemma}
\begin{proof} Write ${\b w}=r{\b v}_1+ s{\b v}_2$ with $r,s\in A$. Then
\[
\begin{aligned}
\deg B(\b w,{\b v}_{3})& \le \sup\{\deg r+\deg e, \deg s+\deg f\}\\
&<\sup\{2\deg r+\mu_{1},2\deg s+\mu_{2}\}\\
&=\deg Q(\b w).
\end{aligned}
\]
\end{proof}

Our next task is to show that by modifying suitably the reduced
bases, the last columns of the matrices in \eqref{E:commonbinary}
can be made equal.

\subsubsection*{The case $\mu_{1}<\mu_{2}<\mu_{3}$, $\mu_1\equiv \mu_2\pmod 2$}

\begin{theorem}\label{T:case1} Let $(L,Q)$ and $(L',Q')$ be isospectral ternary lattices
with strictly increasing minima sequence $\mu_1<\mu_2<\mu_3$
and $\mu_1\equiv \mu_2\pmod 2$. Then they are equivalent.
\end{theorem}
\begin{proof} Let $S$ and $S'$ be their corresponding Gram matrices
as in \eqref{E:commonbinary}. 
Let $M=A{\b v}_1+ A {\b v}_2\subset L$. Since $Q$ represents
$g'$, there exists
${\b v}  \in L$ such that $Q({\b v} )=g'$. Note that for parity reasons
${\b v} \not\in M$, so it is of the form ${\b v} = \lambda {\b v}_{3} + \b w$ 
with $\lambda \in \h^{\times}$ and $\b w\in M$ with $\deg Q(\b w)<\mu_{3}$.
We have
\[
g'=Q({\b v} )=\lambda^{2} g+\lambda B(\b w,{\b v}_{3})+Q(\b w).
\]
Comparing leading coefficients we have $\lambda^{2}=1$. Hence
\[
g'-g= Q(\b w)\pm B(\b w,{\b v}_{3}).
\]
If $\b w\ne 0$, then by Lemma \ref{L:ineq} we get $\deg(g'-g)=\deg Q(\b w)\ge \mu_{1}$,
which contradicts the inequality $\deg(g'-g)< \deg(e'-e)<\mu_{1}$ of Lemma \ref{L:adjoints}.

Thus $g=g'$ and $e^{2}={e'}^{2}$. If $e=e'$, then 
$f=f'$ by \eqref{E:adj} and we are done. So assume $e'=-e\ne 0$ and fix $z\in \h^\times$
such that $b+z e\ne 0$ (the reason for this choice of $z$ will become apparent below). 
Since $Q$ and $Q'$ represent 
in particular the same polynomials, for each $x\in \h$ the equation
\begin{equation}\label{E:eqn2}
Q(x {\b v}_1+{\b v}_2+ z{\b v}_3)=Q'(u {\b v}_1+ v{\b v}_2+z{\b v}_1)
\end{equation}
has a polynomial solution $(u,v)$. Subtracting $z^2 g$ from both sides and using Lemma \ref{L:ineq} we conclude that $\deg Q'(u {\b v}_1+ v{\b v}_2)=\mu_2$, so $v\in \h$. 

Suppose first that for some $x\in \h$, there is a solution $(u,v)$ to \eqref{E:eqn2} with
$v=1$.  Then we have
\begin{equation}\label{E:eqn3}
(x^2-u^2)a +2 (x-u)b +2(f-f')+2  e (x+u) z=0
\end{equation}
Since $\deg (f-f')<\deg e$ by Lemma \ref{L:adjoints}, the above equality implies $x^2=u^2$. If $x=u$,
then \eqref{E:eqn3} reduces to
\[
f - f'=-2 e u
\]
and for degree reasons we must have $f=f'$. By \eqref{E:adj} we get $b=0$
and the transformation ${\b v_2}\mapsto -{\b v_2}$ takes $S$ into $S'$. 
If $x=-u$, then \eqref{E:eqn3} reduces to
\[
f -f'= 2 b u
\]
which similarly implies $f=f'$ since $\deg(f -f')<\deg b$ by \eqref{E:adj}. We conclude as 
in the previous case.

Assume now that for all $x$ all solutions $(u,v)$ to \eqref{E:eqn2} have $v\ne 1$. Then, by
the pigeonhole principle, there must be a pair $(x_1,x_2)\in \h^2$, $x_1\ne x_2$, 
such that the equations
\begin{equation}\label{E:eqn4}
\begin{cases}
Q(x_1 {\b v}_1+{\b v}_2+ z{\b v}_3) = & Q'(u_1 {\b v}_1+ v{\b v}_2+z{\b v}_1)\\
Q(x_2 {\b v}_1+{\b v}_2+ z{\b v}_3) = & Q'(u_2 {\b v}_1+ v{\b v}_2+z{\b v}_1)
\end{cases}
\end{equation}
have solutions $(u_1,v)$ and $(u_2,v)$ (with a common $v$).
Taking the difference of the two equations in \eqref{E:eqn4}, we get
\begin{equation}\label{E:eqn5}
(x_1 - x_2) [a (x_1 + x_2) + 2 b + 2 e z]= (u_1 - u_2) [a (u_1 + u_2) + 2 b v - 2 e z],
\end{equation}
and comparing degrees we see that $u^2_1-u_2^2=x^2_1-x_2^2$. In particular $u_1$ and $u_2$
must be constant. Taking $u_1=u\in \h$ in  \eqref{E:eqn2} we get
\[
(v^2-1)c =a (x^2-u^2)+ 2 b (x-u v) + 2(f - f' v) z + 2e (u + 2) z.
\]
Since all the terms on the right-hand side have degree $<\mu_2=\deg c$, we must have $v^2=1$. 
Having already excluded the case $v=1$, the only allowed value is $v=-1$. Substituting this value
in \eqref{E:eqn5}, cancelling the equal terms and bringing all terms to one side of the equation,
we get
\[
(u_1 - u_2 + x_1 - x_2) (b + e z)=0.
\]
Since we have taken the precaution of choosing $z\in \h^\times$ so that 
$b + e z\ne 0$, we conclude $u_1 - u_2=-x_1 +x_2$, which combined with the
previously established equality $u^2_1-u_2^2=x^2_1-x_2^2$ yields $u_1=-x_1$ and $u_2=-x_2$. Substituting
in the first equation of $\eqref{E:eqn4}$ we get $f =-f'$. Then the transformation 
$\b v_3\mapsto -\b v_3$ takes $S$ into $S'$.

\end{proof}

\subsubsection*{The case $\mu_{1}=\mu_{2}<\mu_{3}$}

Assume that $(L,Q)$ and $(L',Q')$ are isospectral with $\mu_{1}=\mu_{2}$ and let $S$ and $S'$
be their Gram matrices as in Corollary \ref{C: commonbinary}. 

The first step is to show after a suitable change of basis we can also assume $g=g'$.
 Since $Q$ represents $g'$, 
we can write $Q(r {\b v}_1 + s {\b v}_2 + t {\b v}_3)=g'$. Comparing leading coefficients, we see
$t^2=1$, so there is no loss of generality in assuming $t=1$. Consider the transformation
\[
U=\left [\begin{matrix}  1 & 0 & r\\ 0 & 1 & s\\ 0 & 0 & 1 \end{matrix}\right ] .
\]
Then the matrix $S''=U S U^t$ has the form
\[ 
\left [\begin{matrix}  a & b & E \\ b & c& F \\ E & F & g' \end{matrix}\right ] .
\]
Since $\det(S'')=\det(S')$, we have
\[
Q_0(-F,E)=Q_0(-f',e'),
\]
where $Q_0(X,Y)=aX^2+2b XY+ c Z^2$. Since $Q_0$ is definite, $\deg Q_0(-F,E)= \max\{ 2\deg E+\mu_1, 2\deg F+\mu_1\}$
and since $S'$ is reduced, $\deg Q_0(-f',e')<3 \mu_1$ and hence $\deg E<\mu_1$ and $\deg F<\mu_1$. This shows
that $S''$ is reduced, i.e., we can assume henceforth $g=g'$ without loss of generality.

For each $(x,y,z)\in \h ^3$, the equation 
\begin{equation}\label{E:samerepr}
Q(x {\b v}_1 + y {\b v}_2 + z {\b v}_3)= Q'(x' {\b v}_1 + y' {\b v}_2 + z' {\b v}_3)
\end{equation}
has a solution $(x',y',z')$, where $(x',y',z')$ are {\em a priori} polynomials. By
taking leading coefficients, we see $z^2=z'^2$, so $z'$ is in $\h $. Subtracting
the term $z^2 g= z'^2 g'$ on both sides of \eqref{E:samerepr}, and applying Lemma
\ref{L:ineq}, we get
\[
\deg Q(x {\b v}_1 + y {\b v}_2 )= \deg Q'(x' {\b v}_1 + y' {\b v}_2),
\]
which immediately implies $x',y'\in \h $. 

\begin{lemma}\label{L:span} Assume that $Q$ and $Q'$ are ternary definite
isospectral quadratic forms with $\mu_1=\mu_2$, Gram matrices
as in Corollary \ref{C: commonbinary} and the additional condition
$g=g'$. Then $\spa\{e,f\}=\spa\{e',f'\}$.
\end{lemma}

\begin{proof} We shall show $\spa\{e,f\}\subset \spa\{e',f'\}$ and conclude by symmetry. 
If $e=f=0$ there is nothing to prove, so assume $(e,f)\ne (0,0)$. Fix $(x,y)\in \h ^2$ such
that $x e+y f\ne 0$. Then, for all $z\in \h $, the equation
\begin{equation}\label{E:basiceqn}
Q(x,y,z)=Q'(u,v,z)
\end{equation}
has a solution $(u,v)\in \h ^2$. Taking leading coefficients, we see that $(u,v)$ must
satisfy $u^2-\delta v^2=x^2-\delta y^2$, that is, there are at most $q+1$ possible 
pairs $(u,v)$. Up to sign, there are at most $(q+1)/2$ possibilities for $(u,v)$. Since
$q>(q+1)/2$, and the left-hand side of \eqref{E:basiceqn} takes $q$ different values
as $z$ runs over $\h $, there must be $z_1\ne z_2$ in $\h $ and $(u,v)\in \h ^2$
such that
\begin{equation}\label{E:basiceqn2}
Q(x,y,z_1)=Q'(u,v,z_1)\quad \mathrm{and}\quad Q(x,y,z_2)=Q'(\epsilon u,\epsilon v,z_2), 
\end{equation}
where $\epsilon=\pm 1$. Subtracting the two equations we get
\begin{equation}\label{E:linrel}
(z_1-z_2)( xe+y f)=(z_1-\epsilon z_2)(u e' +v f'),
\end{equation}
which shows $ xe+y f\in\spa\{e',f'\}$.

\end{proof}

\begin{lemma}\label{L:auto} Let $Q_0$ be a binary definite quadratic form with $\mu_1=\mu_2$
and let $a$ be a polynomial of degree $\mu_1$ represented by $Q_0$. Then $\Aut(Q_0)$
acts transitively on the set $\{(x,y)\in \h ^2 : Q_0(x,y)=a\}$.
\end{lemma}

\begin{proof} We can assume $Q_0= a X^2+ 2b XY + c Y^2$, where $a,b,c$ are relatively prime 
and $\deg(a)=\deg(c)>\deg(b)$. 

If $a, b, c$ are linearly
independent over $\h $, then $Q_0(x,y)=a$ implies $x=\pm 1$ and $y=0$. We get a similar
conclusion if  $b=0$ and $a,c$ are linearly independent. If $b=0$ and $a$ is proportional to $c$,
$Q_0$ is a multiple of a form over $\h $ and the result is well known. The only case left is 
when $a,b,c$ are linearly dependent and $b\ne 0$. In this case, write $c=-\delta a - 2 m b$,
where $\delta$ is a non-square and $m\in \h $, $m\ne 0$. Suppose $Q_0(x,y)=a$. If $y=0$ we are done, so
we may assume $y\ne 0$. We have $Q_0(x,y)=a(x^2-\delta y^2)+ 2b y(x-m y)$, so by
linear independence of $a$ and $b$ we get the relations
\[
x^2-\delta y^2=1 \quad\mathrm{and} \quad x-m y=0,
\]
which ensure that $U=\left [\begin{matrix}  x & \delta y\\ y & x \end{matrix}\right ] $ is an
automorphism of $Q_0$.
\end{proof}

\begin{lemma}\label{L:hasse}
Let $F, G\in \h [X]$ be  polynomials of degree $2$ 
such that $F(x)\equiv  G(x) \pmod {{\h ^*}^2}$ for all $x\in \h $. Then
$F= u^2 G$, where $u\in \h $. 
\end{lemma}
\begin{proof} The hypothesis implies in particular that
the polynomials $F$ and 
$G$ have the same roots in $\h $ (if any) so there is no loss of generality in assuming
that they are irreducible.

If $F$ and $G$ are relatively prime, then the equation 
$Y^2= F(X) G(X)$ defines an elliptic curve with at least $2 q$ points
over $\h $. This contradicts  Hasse's bound \cite{Rosen:2002kv}[Chapter V] if $q>5$.
For $q=3,5$ the assertion is easily verified by direct computation.
\end{proof}

\begin{theorem}\label{T:case2} If two $Q$ and $Q'$ ternary definite quadratic forms are
isospectral  with $\mu_1=\mu_2$ or $\mu_2=\mu_3$,  then they are equivalent.
\end{theorem}
\begin{proof} If $\mu_2=\mu_3$, we replace $(L,Q)$ and $(L',Q')$ by their adjoints 
which in this case satisfy $\mu_1(Q^\ad)=\mu_1+\mu_2= \mu_1+\mu_3=\mu_2(Q'^\ad)$. So we
can limit ourselves to the case $\mu_1=\mu_2$. 

With the notation and hypotheses of Lemma \ref{L:span},
let $E=\spa\{e,f\}=\spa\{e',f'\}$. If $\dim E=0$, there is nothing to prove;
we will deal with the other two cases.

Suppose first $\dim E=1$. Let $h\in E$ be a monic polynomial of degree $d$
and write $e=e_dh$,
$f=f_d h$, $e=e_d'h$, $f=f_d' h$, where $e_d, f_d, e_d', f_d',$ are in $\h $. 

Let $Q_0(X,Y)=Q(X,Y,0)=Q'(X,Y,0)$. The equality $\det(Q)=\det(Q')$ implies $Q_0(-f,e)=Q_0,-f',e')$.
Dividing by $h^2$ we get $Q_0(-f_d,e_d)=Q_0(-f'_d,e'_d)$. Applying
Lemma \ref{L:auto}, there is an automorphism $U$ of $Q_0$ such that 
$U\left [\begin{matrix} -f \\ e \end{matrix}\right ] = 
\left [\begin{matrix}  -f'\\ e' \end{matrix}\right ] $. Then
\[
W=\left [\begin{matrix}  U & 0 \\ 0 &1 \end{matrix}\right ] 
\]
satisfies $Q W = Q'$, as desired. \\

Suppose now $\dim E=2$. By Lemma \ref{L:span}, there exists a matrix $M\in \GL_2(\h )$
such that 
\[
\left [\begin{matrix}  -f\\ e \end{matrix}\right ] = M \left [\begin{matrix}  -f' \\ e' \end{matrix}\right ] .
\]
We shall prove that $M$ is an automorphism of $Q_0$. Let $(x,y)\in \h ^2$, $(x,y)\ne (0,0)$, and let $(u,v) \in \h ^2$ 
and $z_1, z_2$ as in the proof of Lemma \ref{L:span}. We get from \eqref{E:linrel}

\[
\left [\begin{matrix}  u \\ v \end{matrix}\right ] = h_{x,y} M \left [\begin{matrix}   x \\ y  \end{matrix}\right ] 
\]
with $h_{x,y}=(z_1-z_2)/(z_1-\epsilon z_2)\in \h ^\times$ (depending {\em a priori} on $(x,y)$). 

Let $R=X^2-\delta Y^2$. Since $R(x,y)=R(u,v)$, substituting we have $R(x,y)= h_{x,y}^2 (R M) (x,y)$.
Thus the quadratic forms $R$ and $R M$ represent the same elements of $\h $ up to squares, i.e.
the quadratic polynomials $F(t)=R(t,1)$ and $G(t)=(R M)(t,1)$ satisfy the hypothesis of
Lemma \ref{L:hasse}, hence  $R=s^2 R M$ for some $s\in \h ^\times$, i.e. $h_{x,y}^2=s^2$ for all $(x,y)\in \h ^2$.

Now from the equality $\det(Q)=\det(Q')$, we get $Q_0 (-f,e)= Q_0 (-f',e')$. Let $d=\max\{\deg e,\deg f\}$ and
take coefficients of degree $\mu_1 +2 d$ in this equality. Then
\[
R(-f_d,e_d) = R(-f_d',e_d') \ne 0,
\]
and therefore $s^2=1$ and $R= R M$.

If $h_{x,y}=1$, we conclude from the first equation in 
\eqref{E:basiceqn2} that $Q_0(x,y)=Q_0(u,v)$. If $h_{x,y}=-1$,
then $\epsilon =-1$ and $z_1=0$ and we conclude again from
\eqref{E:basiceqn2} that $Q_0(x,y)=Q_0(u,v)$. Thus $Q_0 = Q_0 M$;
this condition ensures that
\[
N:=\left [\begin{matrix}  M & 0 \\ 0 &1 \end{matrix}\right ] 
\]
satisfies $Q'=Q N$.

\end{proof}

\subsubsection*{The case $\mu_{1}<\mu_{2}<\mu_{3}$,  $\mu_1\equiv \mu_3 \pmod 2$} 

Let $(W,\phi)$ be a quadratic space over $\h $ of dimension $n$ and rank $r$. Recall that the Gauss sum associated
to $(W,\phi)$ is defined by
\[
\Gamma(W,\phi)=\sum_{w\in W} \chi(\phi(w)),
\]
where $\chi: \h  \to \C^\times$ is the character defined by $\chi(u)=\exp(2\pi i\Tr(u)/p)$ and
$\Tr:\F_q\to \F_p$ is the trace to the prime field $\F_p$.

It is immediate from the definition that $\Gamma$ is multiplicative on orthogonal sums.  
Let $W_1=\rad(W,\phi)$ and let $W_0\subset W$ be a complement of $W_1$. Then
$\Gamma(W,\phi)=\Gamma(W_0,\phi_0)\Gamma(W_1,0)$, where $\phi_0=\phi|_{W_0}$. Clearly 
$\Gamma(W_1,0)=q^{n-r}$. Writing $\phi_0=\sum_{i=1}^r a_i X_i^2$ in some orthogonal 
basis of $W_0$, we get $\Gamma(W_0,\phi_0)=\Gamma(\F_q,\langle a_1\rangle )\cdots \Gamma(\F_q,\langle a_r\rangle )$.
Using further the property that $\Gamma(\F_q,\langle a_i\rangle )=\psi(a_i) G$, where $G= \Gamma(\F_q,\langle 1\rangle )$ 
and $\psi:\h ^\times \to \{\pm 1\}$ is the quadratic character (see e.g. \cite[Proposition 6.3.1]{Ireland:1982lr}), we get
\begin{equation}\label{E:Gamma}
\Gamma(W,\phi)=
q ^{n-r} \psi (\det \phi_0) G^r. 
\end{equation}

Note that in particular, $\Gamma(W,\phi)=\Gamma(W,\phi')$ if and only if $(W,\phi)\simeq (W,\phi')$.

\begin{definition} Let $\Phi=(\phi_1,\phi_2,\cdots,\phi_m)$ and $\Phi'=(\phi'_1,\phi'_2,\cdots,\phi'_m)$ be systems of quadratic forms on $W$, i.e. quadratic mappings $W\to \h ^m$. We shall say that $\Phi$ and $\Phi'$ are {\em isospectral} if
$|\Phi^{-1} (\b y)|= |{\Phi'}^{-1} (\b y)|$ for all $\b y\in \h ^m$.
\end{definition}

The following theorem is a particular case of a result by Carlitz \cite[Theorems 3.2-3.3]{Carlitz:1954tg}
on systems of polynomial equations.

\begin{theorem}[Carlitz] \label{T: Carlitz} Two systems of quadratic forms $\Phi$ and $\Phi'$
as above are isospectral if and only if

\[
\Gamma(\sum_{i=1}^m x_i \phi_i)= \Gamma(\sum_{i=1}^m x_i \phi_i')
\]
for all $(x_1,x_2,\ldots,x_m)\in \h ^m$.

\end{theorem}

Let $Q, Q'$ be isospectral definite quadratic forms with successive minima $(\mu_1,\mu_2,\mu_3)$ on $L$ and let $W=L_{\mu_3}$. Write $Q(\b x)=\sum_{i=0}^{\mu_3}
Q_i(\b x) t^i$ (respectively $Q'(\b x)=\sum_{i=0}^{\mu_3}
Q_i'(\b x) t^i$). Then the systems $\Phi=(Q_0,\ldots,Q_{\mu_3})$
and $\Phi'=(Q'_0,\ldots,Q'_{\mu_3})$ are isospectral. Let $B$, $B_i$, $B'$, $B_i'$ be
the symmetric bilinear forms associated to $Q$, $Q_i$, $Q'$, $Q_i'$. By Theorem
\ref{T: Carlitz} and \eqref{E:Gamma}, we have in particular
\begin{equation}\label{E:dets}
\det (\sum_{i=0}^{\mu_3} x_i B_i) \equiv \det (\sum_{i=0}^{\mu_3} x_i B'_i) \pmod {{\h ^\times}^2}
\end{equation}
for all $(x_0,x_2,\ldots,x_{\mu_3})\in \h ^{\mu_3 +1}$.

Let $k_1=(\mu_3-\mu_1)/2$ and $k_2= (\mu_3-\mu_2-1)/2$. We fix the basis
\begin{equation}\label{E:basis}
\{{{\b v}_1}, t{{\b v}_1}, \ldots, t^{k_1} {{\b v}_1}, {{\b v}_2},t {{\b v}_2},\ldots, t^{k_2} {{\b v}_2}, {{\b v}_3}\}
\end{equation}
 of $W$ and identify all the symmetric bilinear forms on $W$ with their respective matrices in this basis.

\begin{lemma}\label{L:dets} With the notation above, we have

\[
\det (\sum_{i=0}^{\mu_3-1} x_i B_i) = \det (\sum_{i=0}^{\mu_3-1} x_i B'_i)
\]
for all $(x_0,x_2,\ldots,x_{\mu_3-1})\in \h ^{\mu_3}$.
\end{lemma}
\begin{proof}

Fix $(x_0,x_2,\ldots,x_{\mu_3-1})\in \h ^{\mu_3}$ and consider $\det (\sum_{i=0}^{\mu_3} x_i B_i)$ 
and $\det (\sum_{i=0}^{\mu_3} x_i B'_i)$ as polynomials in the variable $x_{\mu_3}$. They
have degree two in $x_{\mu_3}$, the same leading coefficient ($= -\delta$) and are equal
up to squares of $\h ^\times$ by \eqref{E:dets}, so, by Lemma \ref{L:hasse}, they must be equal  as polynomials in $x_{\mu_3}$. We conclude by taking
$x_{\mu_3}=0$. 
\end{proof}

\begin{lemma} Let $m=(\mu_1+\mu_3)/2$. Then for all $m\le j\le \mu_3$ we have $B'_j=B_j$.
\end{lemma}
\begin{proof} For $\b x= (x,y,z)\in W$, we have $Q(\b x)-Q'(\b x)= 2(e-e') xz+ 2(f-f') yz + (g-g') z^2$.
By Lemma \ref{L:adjoints}, all three terms have degrees $<m$.
\end{proof}

We shall use the following notation henceforth: $n=\dim W$, $s=\max\{m-\mu_2,-1\}$,  $r=k_2-s$. (Note
that $n=k_1+k_2+3=(k_1+1)+(s+1)+r+1$.) 

\begin{lemma}\label{L:properties} The forms $B_j$ have the following properties
\begin{enumerate}
\item $B_l(t^i{\b v}_1, t^j{\b v}_1)=0$ for $l\ge m$ and $i+j<k_1$
\item $B_m(t^i{\b v}_1, t^j{\b v}_1)=1$ for $i+j=k_1$
\item $B_l(t^i{\b v}_1, t^j{\b v}_2)=0$ for $l-i-j\ge \mu_1$
\item $B_l(t^i{\b v}_2, t^j{\b v}_2)=0$ for $l\ge m$ and $i+j<s$
\item $B_m(t^i{\b v}_2, t^j{\b v}_2)=c_{\mu_2}$ for $i+j=s$
\item $B_l(t^i{\b v}_1, {\b v}_3)=0$ for $l\ge m$ and $i\le k_1$
\item $B_l(t^i{\b v}_2, {\b v}_3)=0$ for $l\ge m$ and $i\le s$.
\end{enumerate}
\end{lemma}
\begin{proof} The lemma follows immediately from the fact that $\{{\b v}_1,{\b v}_2, {\b v}_3\}$ is a reduced 
basis for $Q(\b x)$. 
\end{proof}

Let
$
{\mathcal B}=X_m B_m + \sum_{j=0}^{r-1}  X_{\mu_3-1-2j} B_{\mu_3-1-2j}.
$
It follows from Lemma \ref{L:properties} that the matrix of ${\mathcal B}$ in the basis \eqref{E:basis} has the form
\[
\label{E:matrix1}
\includegraphics[scale=0.65]{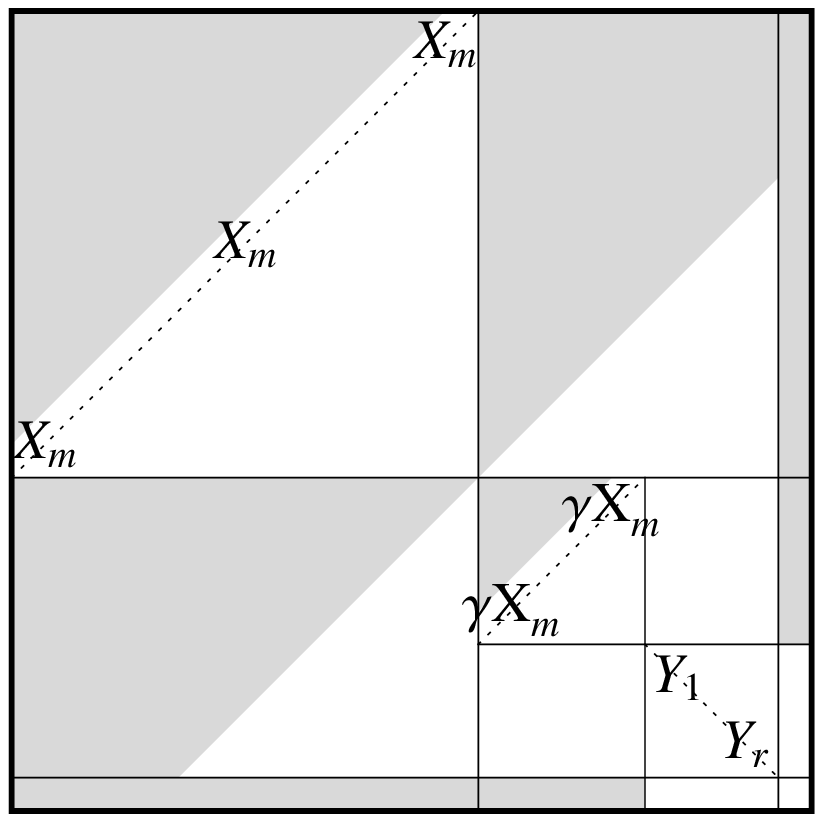}
\]
where the greyed areas consist entirely of zeros and the sizes of the blocs correspond to the
partition $n=(k_1+1)+(s+1)+(r)+(1)$. Here $Y_j= {\mathcal B}(t^{k_2+j-r} {{\b v} }_2,t^{k_2+j-r} {{\b v} }_2)$ for $j=1,\ldots,r$
and $\gamma=c_{\mu_2}$.

\begin{lemma}\label{L:minors} Let ${\mathcal C} =(\rho_{ij})$ be the adjoint of the matrix ${\mathcal B}$
and let $M= \break X_m^{k_1+k_2+2-r} \prod_{j=0}^{r-1}  X_{\mu_3-1-2j}$.
\begin{enumerate}
\item 
When $i<n$, the coefficient of $M$ in $\rho_{ni}$ is equal to $0$
\item 
The coefficient of $M$ in $\rho_{nn}$ is equal to $\pm \gamma^{k_2+1}$.
\end{enumerate}
(Note that the entries of ${\mathcal C}$ are homogeneous polynomials of degree $n-1$ in the variables $X_j$.)
\end{lemma}
\begin{proof}
Leaving the bottom row untouched, we apply elementary row operations 
to clear the entries below the ``diagonals'' containing $X_m$. This can be accomplished
in the ring $\h [X_m^{\pm 1}][X_{\mu_3-2r+1},\cdots, X_{\mu_3-1}]$. We get a matrix of the form

\[
\label{E:matrix2}
\includegraphics[scale=0.65]{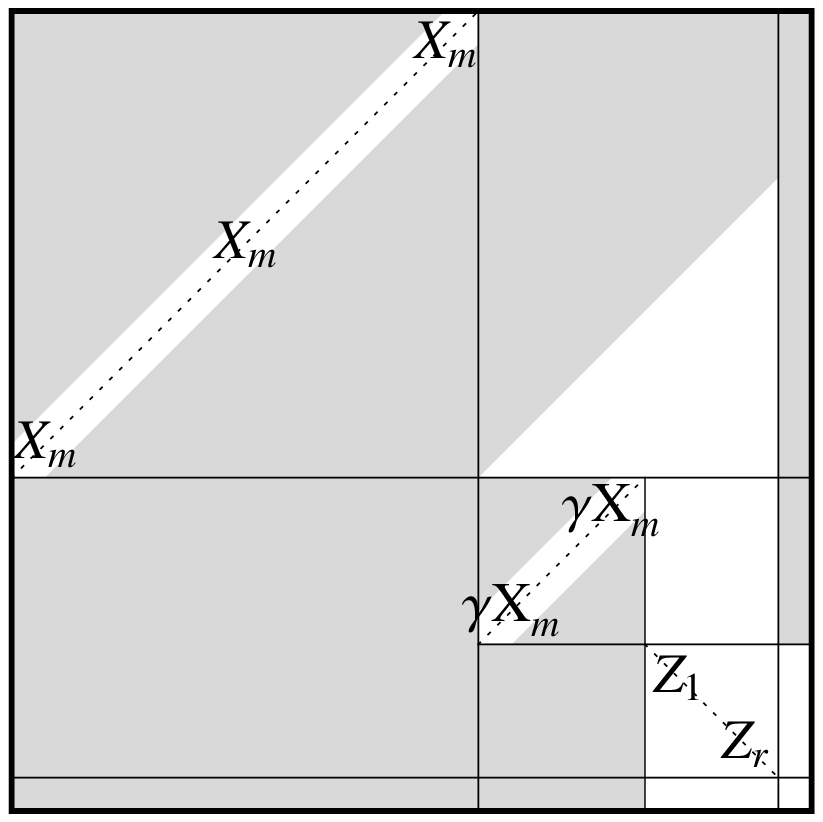},
\]

Note that all the coefficients above (and to the left of ) $Y_j$ are linear combinations of variables $X_i$ with 
$i< 2j+2m-\mu_2$, so $Z_j$ and $Y_j$ have the same term in $X_{2j+2m-\mu_2}$, namely $\gamma X_{2j+2m-\mu_2}$.

The elementary row operations have not altered the minors of ${\mathcal B}$ along the bottom row (i.e. the determinants
of the submatrices obtained by removing the bottom row and a column). It is clear from the shape of the above matrix that in these minors, only the product
$Y_1\cdots Y_r$ can yield a term divisible by  $\prod_{j=0}^{r-1}  X_{\mu_3-1-2j} $. Thus the minors $\rho_{ni}$ obtained by
removing a column different from the last one ($i<n$) do not contain monomials divisible by $\prod_{j=0}^{r-1}  X_{\mu_3-1-2j} $.
The coefficient of $M$ in the minor $\rho_{nn}$ is $\pm \gamma^{r+s+1}$.

\end{proof}

\begin{lemma}
For all $0\le i<m$ we have
\[
\det( X_i B_i' +{\mathcal B})-\det (X_i B_i +{\mathcal B}) =\pm (g_i'-g_i) \gamma^{r+s+1} X_i M + N,
\]
where $M=X_m^{k_1+k_2+2-r} \prod_{j=0}^{r-1}  X_{\mu_3-1-2j} $ and $N$ is divisible by 
 $X_i^2$.

\end{lemma}

\begin{proof}
Expanding as polynomials in $X_i$ and separating the linear part, we have

\begin{equation}\label{E:linearpart}
\begin{aligned}
\det( X_i B_i' +{\mathcal B})-\det (X_i B_i +{\mathcal B})=&\Tr({\mathcal C} (B_i'-B_i)) X_i \\
& + \mathrm{terms\ divisible\  by\ } X_i^2.
\end{aligned}
\end{equation} 
The matrix of $B_i'-B_i$ with respect to the basis \eqref{E:basis} has zeros everywhere except possibly on the last row and the last column and $(B_i'-B_i)_{nn}= g_i'-g_i$. Combining this with Lemma \ref{L:minors} we get
that the coefficient of $M$ in $\Tr({\mathcal C} (B_i'-B_i))$ is $\pm \gamma^{r+s+1} (g_i'-g_i)$. The lemma follows now
immediately from \eqref{E:linearpart}.

\end{proof}

\begin{corollary}\label{C:fund} If $q> k_1+k_2+2-r$, then $g=g'$.

\end{corollary} 

\begin{proof} A monomial  of $N$ that is equal to $X_j M$ as functions on $\h $ must
be of the form $X_j^{q^s} P$, where $P$ is divisible by {\em all} the variables
other than $X_j$. In particular $\deg P\ge r+1$, so $q^s\le \dim L_{\mu_3}-(r+1)=k_1+k_2+2-r$,
which implies $s=0$. Since $\det( X_i B_i' +{\mathcal B})=\det (X_i B_i +{\mathcal B})$ as functions, we must have
$g_i=g'_i$ for $0\le i<m$. Since $\deg(g-g')< \mu_1<m$ by Lemma \ref{L:adjoints}, we must have $g=g'$.

\end{proof} 

\begin{theorem}\label{T:case3} If $q> \max\{2+\mu_3-\mu_2, 2+\mu_2-\mu_1\}$ then $Q$ and $Q'$ are isometric.

\end{theorem}

\begin{proof} The condition on $q$ ensures that both pairs $(Q, Q')$ and $(Q^\ad, {Q'}^\ad)$
satisfy the hypotheses of Corollary \ref{C:fund}. Applying  Corollary \ref{C:fund} to 
$(Q, Q')$ we get $g=g'$ and hence $e^2=e'^2$. Applying it to  $(Q^\ad, {Q'}^\ad)$, we get 
$cg-f^2=cg'-f'^2$ and hence $f^2=f'^2$. 

There is no loss of generality in assuming $e=e'$.
If $f=f'$ we are done, so assume $f=-f'\ne 0$. Comparing determinants we get $be=0$. If $b=0$,
then the transformation ${\b v}_2\mapsto -{\b v}_2$ changes $f$ into $-f$ and leaves the rest alone. Similarly,
if $e=0$, the transformation ${\b v}_3\mapsto -{\b v}_3$ changes $f$ into $-f$ and leaves the other coefficients
unaltered.

\end{proof} 

\nocite{Bureau:2006fk}

\bibliographystyle{amsplain}

\providecommand{\bysame}{\leavevmode\hbox to3em{\hrulefill}\thinspace}

\providecommand{\href}[2]{#2}

\end{document}